\documentclass[9pt,amsfonts]{article}
\usepackage{graphicx,amssymb,eucal,mathrsfs}
\usepackage{latexsym,amsmath,epsfig,epic,eepic,xcolor}
\usepackage[all]{xy}
\usepackage{wasysym,ntheorem,hyperref}
\usepackage{epsf,graphicx,pgf,etex,amscd,pstricks}

\newcommand{\Z}{\mathbb{Z}}
\newcommand{\R}{\mathbb{R}}

\newcommand{\C}{\mathbb{C}}

\newcommand{\PSLR}{\mathrm{PSL}_{2}(\mathbb{R})}

\newtheorem{theorem}{{\bf Theorem}}

\newtheorem{prop}[theorem]{Proposition}
\newtheorem{lemma}[theorem]{Lemma}
\newtheorem{rem}{Remark}
\newtheorem{cor}{Corollary}

\theoremnumbering{Alph}
\newtheorem{main}{{\bf  Theorem}}

\title{Coxeter groups and K\"{a}hler groups}
\author{Pierre Py}
\date{June 2013}

\begin{document}

\maketitle

\begin{abstract} We study homomorphisms from K\"{a}hler groups to Coxeter groups. As an application, we prove that a cocompact complex hyperbolic lattice (in complex dimension at least $2$) does not embed into a Coxeter group or a right-angled Artin group. This is in contrast with the case of {\it real} hyperbolic lattices. 
\end{abstract}


\section{Introduction}

A {\it K\"ahler group} is by definition the fundamental group of a compact K\"{a}hler manifold. We refer the reader to~\cite{abckt} for an introduction to these groups, and to~\cite{biswasmj,burger,kotschick} for more recent developments. The purpose of this note is to study homomorphisms from K\"ahler groups to Coxeter groups (for the definition of Coxeter groups, see Section~\ref{rappels} and~\cite{davisbook,humphreys} for a more detailed introduction). 

We first recall a few classical definitions. In this text, a $2$-dimensional orbifold $\Sigma$ is a compact Riemann surface $S$ marked with a finite set of points $p_{1},\ldots , p_{n}$, each point $p_{i}$ being assigned a multiplicity $m_{i}\ge 2$. The orbifold $\Sigma$ is called {\it hyperbolic} if its orbifold Euler characteristic 
$$\chi^{orb}(\Sigma):=\chi(S)-\sum_{i=1}^{n} \left(1-\frac{1}{m_{i}}\right)$$ 

\noindent is negative. In this case, there exists a cocompact discrete subgroup $\Gamma \subset \PSLR$ acting on the unit disc $\Delta \subset \C$, such that $S$ can be identified with the quotient $\Delta/\Gamma$, the $p_{i}$ corresponding to orbits of points of $\Delta$ with a non-trivial stabilizer in $\Gamma$. The group $\Gamma$ is isomorphic to the orbifold fundamental group $\pi_{1}^{orb}(\Sigma)$ of $\Sigma$, i.e. the quotient of the group $\pi_{1}(S-\{p_{1}, \ldots , p_{n}\})$ by the normal subgroup generated by the elements $\gamma_{i}^{m_{i}}$ ($\gamma_{i}$ being the conjugacy class of a loop around the puncture $p_{i}$). In such a situation, we will always think of $\Sigma$ as a quotient of the unit disc and say that $\Sigma$ is a {\rm hyperbolic $2$-orbifold}. 

A map from a complex manifold $X$ to a hyperbolic $2$-orbifold $\Sigma$ is holomorphic if it lifts to a holomorphic map from the universal cover of $X$ to the unit disc. A fibration from $X$ onto $\Sigma$ is a proper holomorphic surjective map $f : X \to \Sigma$ with connected fibers; such a map induces a surjection $f_{\ast} : \pi_{1}(X)\to \pi_{1}^{orb}(\Sigma)$ (see for instance Lemma 3 in~\cite{catanese} for more on this topic). 

In the following, if $G$ is any group, we will denote by $G_{{\rm ab}}$ the abelianization of $G$. If $f_{1} : G\to H_{1}$ and $f_{2} : G\to H_{2}$ are homomorphisms we will say that {\it $f_{1}$ factors through $f_{2}$} if the kernel of $f_{2}$ is contained in the kernel of $f_{1}$.

\newpage

We can now state the:

\begin{main}\label{kc} Let $X$ be a compact K\"ahler manifold with fundamental group $\Gamma=\pi_{1}(X)$. Let $W$ be a Coxeter group and $\varphi : \Gamma \to W$ be any homomorphism. Then, there is a finite cover $$X_{0}\to X$$ with fundamental group $\Gamma_{0}:=\pi_{1}(X_{0})$, and finitely many fibrations $p_{i} : X_{0}\to \Sigma_{i}$ ($1\le i\le N$) onto hyperbolic $2$-orbifolds such that the restriction of $\varphi$ to $\Gamma_{0}$ factors through the map 
$$\Gamma_{0}\to (\Gamma_{0})_{{\rm ab}} \times \pi_{1}^{orb}(\Sigma_{1})\times \cdots \times \pi_{1}^{orb}(\Sigma_{N})$$
induced by the $p_{i}$'s and by the natural map $\Gamma_{0}\to (\Gamma_{0})_{{\rm ab}}$. 
\end{main}

From now on, by a {\it surface group} we will mean the fundamental group of a closed orientable surface of genus greater than $1$. The orbifold fundamental groups appearing above have finite index subgroups isomorphic to surface groups. One thus deduces from Theorem~\ref{kc} that if a K\"ahler group $\Gamma$ admits a {\it faithful} homomorphism into a Coxeter group, it must have a finite index subgroup $\Gamma_{1}$ isomorphic to a subgroup of the direct product of a free Abelian group with a direct product of surface groups. We refer the reader to~\cite{delzantgromov} for other situations where one can construct faithful homomorphisms from certain K\"ahler groups to products of surface groups (possibly with an Abelian factor).

On the other hand, Bridson, Howie, Miller and Short~\cite{bhms2002} have studied subgroups of direct products of surface groups and free groups, from the point of view of their homological properties; see also~\cite{bridsonhowie,bhms2009} for more general results. They proved in~\cite{bhms2002} that a subgroup of a direct product of free and surface groups with strong enough finiteness properties is virtually isomorphic to a product of finitely generated subgroups of the factors (see Section~\ref{proofs} for a precise statement). Hence, it is tempting to use their result to obtain restrictions on K\"ahler groups with strong enough finiteness properties (this possibility was already mentioned in~\cite{bridsonhowie,bhms2009}). Indeed, from the Theorem above and the results of~\cite{bhms2002} we easily deduce (see Section~\ref{proofs}):

\begin{cor}\label{coxcom} If a Coxeter group $W$ is commensurable with a K\"ahler group, then any infinite irreducible factor of $W$ is either Euclidean or has a finite index subgroup isomorphic to a surface group. 
\end{cor}

From Theorem~\ref{kc}, one also deduces immediately:

\begin{cor}\label{noninj} Let ${\rm PU}(n,1)$ be the group of holomorphic isometries of the unit ball in $\C^{n}$ and $\Gamma \subset {\rm PU} (n,1)$ be a cocompact lattice. If $n\ge 2$, $\Gamma$ does not admit any faithful homomorphism into any Coxeter group. 
\end{cor}

Although this last corollary could also be deduced from~\cite{bhms2002}, it admits the following direct proof. First we observe that a subgroup $G$ of a direct product $A\times B$ of torsion-free groups which does not embed in either $A$ or $B$ contains elements of the form $(a,1)$ and $(1,b)$ (for $a, b \neq 1$). In particular, it contains a copy of $\Z^{2}$. Hence a Gromov hyperbolic group embeds in such a direct product only if it embeds in one of the factors. Combined with the remarks above, this proves Corollary~\ref{noninj}.

The statement of Corollary~\ref{noninj} for {\it non-uniform} lattices in ${\rm PU} (n,1)$ follows from the fact that nilpotent subgroups of Coxeter groups are virtually Abelian. Note that higher rank lattices cannot embed in Coxeter groups; this follows from the result in~\cite{noskovvinberg}, combined, for instance, with Margulis' normal subgroup theorem (see~\cite{singh} for a weaker statement).

We now recall the definition of a {\it right-angled Artin group} (denoted by RAAG in what follows). Let $\mathscr{G}$ be a finite graph with vertex set $V(\mathscr{G})$. The RAAG $A(\mathscr{G})$ associated to $\mathscr{G}$ is the group defined by the following presentation: 
$$A(\mathscr{G}):=\big \langle g_{v}, v\in V(\mathscr{G})\; \vert \; [g_{v},g_{w}]=1 \; {\rm if} \; (v,w)\; {\rm is}\; {\rm an}\; {\rm edge}\;  {\rm of}\; \mathscr{G}\big \rangle .$$ 

\noindent See for instance~\cite{charney} for an introduction to RAAG's. Note that Davis and Januszkiewicz~\cite{dj-jpaa} have proved that any RAAG embeds (as a subgroup of finite index) into a right-angled Coxeter group. This implies that Theorem~\ref{kc} and Corollary~\ref{noninj} above still hold when one replaces Coxeter groups by RAAG's. In particular one obtains: 

\begin{cor} Let $\Gamma \subset {\rm PU} (n,1)$ be a cocompact lattice with $n\ge 2$. Then $\Gamma$ does not admit any faithful homomorphism into any right-angled Artin group. 
\end{cor}

This contrasts with the fact that the fundamental group of any compact real hyperbolic $3$-manifold as well as the fundamental groups of ``standard" arithmetic real hyperbolic manifolds in all dimensions  virtually embed into RAAGs. The first statement follows from Agol's recent solution to the virtual Haken conjecture; see~\cite{agol} and the references there. For the case of standard arithmetic hyperbolic manifolds, see~\cite{bhw}. The non-existence of {\it quasi-isometric} embeddings of complex hyperbolic lattices into RAAGs was already known thanks to~\cite{delzantgromov} (see also~\cite{longreid}).

From the previous results we will also deduce:

\begin{cor}\label{corraag} If a {\rm RAAG} $A(\mathscr{G})$ is commensurable with a K\" ahler group, then $A(\mathscr{G})$ is free Abelian of even rank. 
\end{cor}

The fact that the only RAAG's isomorphic to K\"ahler groups are the free Abelian groups of even rank was already established in a different way in~\cite[\S 11.13]{dps2} (that article also describes which RAAG's are fundamental groups of quasi-K\"ahler manifolds).

Note that any RAAG acts properly and cocompactly on a {\it ${\rm CAT}(0)$ cubical complex}, see~\cite{charney} for the definition. Although partial results on actions of K\" ahler groups on ${\rm CAT}(0)$ cubical complexes were obtained by Delzant and Gromov (see~\cite{delzantgromov} or~\cite{burger}), it would be interesting to study these actions in general.

The proof of our results is an easy combination of two classical facts. First, we use the fact that Coxeter groups act faithfully and properly on a product of finitely many trees: this appears for instance in the work of Dranishnikov and Januszkiewicz~\cite{draja} and Januszkiewicz~\cite{janusz2002}, see also~\cite{noskovvinberg}. More recently, this construction was used by L\'ecureux~\cite{lecureux}. Second, we use the fact, due to Gromov and Schoen~\cite{gs}, that a non-elementary action of the fundamental group of a K\"ahler manifold $X$ on a tree gives rise to a fibration of $X$ onto a hyperbolic $2$-orbifold.  

After recalling a few facts about Coxeter groups and their Davis complexes in Section~\ref{rappels}, we prove Theorem~\ref{kc} and Corollaries~\ref{coxcom} and~\ref{corraag} in Section~\ref{proofs}. 

{\bf Acknowledgements.} I would like to thank Jean L\'ecureux who first told me a long time ago about actions of Coxeter groups on products of trees as well as Nicolas Bergeron, Yves de Cornulier, Fr\'ed\'eric Haglund and Luis Paris for useful conversations. I would also like to thank Misha Kapovich; I started thinking about this work after a conversation with him at Oberwolfach in July 2012. Both the question of the embedding of complex hyperbolic lattices into RAAG's and into Coxeter groups are due to him. 

Finally, I would like to thank the referee for his/her detailed reading of the text. 


\section{Coxeter groups and Davis complexes}\label{rappels}

A {\it Coxeter system} $(W,S)$ is a group $W$ with a finite generating set $S\subset W$ all of whose elements have order $2$ and such that $W$ admits the presentation:
$$\langle S\vert (st)^{m_{s, t}}\rangle$$

\noindent where $m_{s, t}\in [1,+\infty]$ is the order of $st$. To the pair $(W,S)$, one associates its {\it Coxeter diagram} $\mathscr{G}$ which is the graph whose set of vertices is parametrized by $S$ and where two vertices $s$ and $t$ are adjacent if $m_{s, t}\ge 3$. We say that $(W,S)$ is irreducible if its Coxeter diagram is connected. In general the partition of $\mathscr{G}$ into connected components gives rise to a decomposition of $S$ as a disjoint union $S = S_{1}\sqcup\cdots \sqcup S_{p},$ and of $W$ as a direct product:

$$W:=W_{1}\times \cdots \times W_{p}$$

\noindent where $W_{i}$ is the subgroup generated by $S_{i}$. Each Coxeter system $(W,S)$ has a canonical faithful linear representation defined as follows. One considers a vector space $V$ with a basis $(u_{s})_{s\in S}$ indexed by the elements of $S$. On $V$ we define a symmetric bilinear form $B$ by:
$$B(u_{s},u_{t})=-\cos \left(\frac{\pi}{m_{s, t}}\right).$$
For each $s\in S$ we define the reflection $\sigma_{s} : V\to V$ by $\sigma_{s}(v)=v-2B(v,u_{s})u_{s}$. There is a unique homomorphism $\sigma : W\to {\rm GL}(V)$ such that $\sigma(s)=\sigma_{s}$; it is faithful and the group $\sigma(W)$ preserves the form $B$. This result is due to Tits, see~\cite[\S 5.3]{humphreys}.

We say that an irreducible Coxeter group $(W,S)$ is Euclidean if the bilinear form $B$ defined above is positive semidefinite but not positive definite. In this case, $W$ can be realized as a cocompact discrete group of affine isometries of a Euclidean space, generated by affine reflections~\cite[\S 6.5]{humphreys}. We finally recall two more facts about Coxeter groups. The first one is the construction of their {\it Davis complex}. The second one concerns decompositions into direct products.

We start with the description of the Davis complex of a Coxeter group. This is a contractible simplicial complex on which the group $W$ acts properly and cocompactly (see~\cite{davisbook} for a detailed study). Say that a subset $T\subset S$ is {\it spherical} if the associated group $W_{T}:=\langle T\rangle \subset W$ is finite. In this case we also say that $W_{T}$ is a spherical special subgroup of $W$. We define $W\mathscr{S}$ to be the union of all cosets of spherical special subgroups:
$$W\mathscr{S} = \bigsqcup_{T} W/W_{T},$$

\noindent where the union runs over the spherical subsets of $S$. The set $W\mathscr{S}$ is partially ordered by inclusion. One gets a simplicial complex ${\rm Flag}(W\mathscr{S})$ (with set of vertices equal to $W\mathscr{S}$) by considering flags in $W\mathscr{S}$: a flag is a finite totally ordered subset of $W\mathscr{S}$, i.e. a finite chain 
$$u_{1}W_{T_{1}}\subset \cdots \subset u_{l}W_{T_{l}}.$$

\noindent The Davis complex $\Sigma(W)$ is the geometric realization of ${\rm Flag}(W\mathscr{S})$. The group $W$ acts naturally by left translation on $W\mathscr{S}$; this action induces a simplicial action of $W$ on $\Sigma(W)$. A {\it reflection} is an element of $W$ which is conjugated to an element of $S$; the {\it wall} associated to a reflection is its fixed point set in $\Sigma(W)$.

The following proposition is classical:

\begin{prop} The complex $\Sigma(W)$ carries a $W$-invariant piecewise Euclidean ${\rm CAT}(0)$ metric. As a consequence, $\Sigma(W)$ is contractible. For any reflection $r\in W$, the fixed point set $Fix(r)\subset \Sigma(W)$ separates $\Sigma(W)$ into two connected components (called the half-spaces associated to $r$).  
\end{prop}
{\it Proof.} The fact that $\Sigma(W)$ carries an invariant ${\rm CAT}(0)$ metric is due to Moussong, see~\cite[Ch. 12]{davisbook}. For a proof of the fact that the fixed-point set of a reflection separates $\Sigma(W)$ into two components, see for instance~\cite[\S 5.3.3]{davisbook}.\hfill $\Box$

Finally, we will need to use the following fact. See~\cite{cornulierharpe,paris,qi} for various proofs of it. 

\begin{prop}\label{nn} Let $W$ be an infinite irreducible Coxeter group. Assume that $W$ is not Euclidean. Let $G\subset W$ be a finite index subgroup. If $G\simeq A\times B$ is decomposed as a direct product of two subgroups $A$ and $B$, either $A=\{1\}$ or $B=\{1\}$.
\end{prop}


\section{Proofs}\label{proofs}

Let $W_{0}$ be a torsion-free finite index normal subgroup of $W$ (such a subgroup exists since Coxeter groups are virtually torsion-free, by Selberg's lemma). We recall here how to produce some actions of $W_{0}$ on certain simplicial trees constructed from the Davis complex $\Sigma(W)$ of $W$. This construction is taken from~\cite{draja,janusz2002}. The key fact is the following observation. 

\begin{center}
{\it If $H$ is a wall of $\Sigma$, and if $\gamma\in W_{0}$, then either $\gamma(H)=H$ or $\gamma(H)\cap H=\emptyset$.}
\end{center}

\noindent See Lemma 1 in~\cite{janusz2002} for a proof. Choose now a $W_{0}$-orbit of walls, say $\mathscr{O}$. We define a tree $T_{\mathscr{O}}$ associated to $\mathscr{O}$ as follows. Let 

$$U=\Sigma(W)-\bigcup_{H\in \mathscr{O}}H.$$ 

\noindent The vertex set of $T_{\mathscr{O}}$ is the set of connected components of $U$; two connected components $U_{1}$ and $U_{2}$ are adjacent if the intersection of their closures is nonempty (in which case it is a wall from $\mathscr{O}$). One obtains in this way a graph. It is easy to see that $T_{\mathscr{O}}$ is a tree: indeed let $e_{H}$ be the edge associated to a wall $H\in \mathscr{O}$. From the fact that the set $\Sigma(W)-H$ has two connected components, one sees that $T_{\mathscr{O}}-e_{H}$ has two connected components. We refer the reader to~\cite{draja,janusz2002} or \cite[\S 14.1]{davisbook} for more details on this construction. Note that $W$ sits inside $W\mathscr{S}$, which is the set of vertices of $\Sigma(W)$. The image of $W$ in $\Sigma (W)$ does not intersect any wall, hence there is a natural map $p_{\mathscr{O}} : W \to T_{\mathscr{O}}$. The group $W_{0}$ acts on $T_{\mathscr{O}}$ and the map $p_{\mathscr{O}}$ is $W_{0}$-equivariant.

Consider now the collection of all $W_{0}$-orbits of walls in $\Sigma(W)$; there are finitely many such orbits $\mathscr{O}_{1}, \ldots, \mathscr{O}_{k}$. Write $T_{i}=T_{\mathscr{O}_{i}}$ and $p_{i}=p_{\mathscr{O}_{i}}$ for the corresponding trees and projections. We get a map 
$$F =(p_{1},\ldots , p_{k}) : W\to T_{1}\times \cdots \times T_{k},$$

\noindent which is proper (see for instance~\cite{janusz2002}). Since $W_{0}$ is torsion-free, the properness of $F$ implies:    

\begin{lemma} The action of $W_{0}$ on $T_{1}\times \cdots \times T_{k}$ is free.
\end{lemma}

We are now ready to prove our main result. 

{\it Proof of Theorem~\ref{kc}.} We consider a homomorphism $\varphi : \Gamma \to W$ where $\Gamma= \pi_{1}(X)$ is K\"ahler and $W$ is a Coxeter group. Let $W_{0}$ be a torsion-free normal subgroup of finite index of $W$ and put $\Gamma_{0}:=\varphi^{-1}(W_{0})$. Let $T_{1}, \ldots , T_{k}$ be the simplicial trees obtained from the construction above. Via the homomorphism $\varphi$, the group $\Gamma_{0}$ acts isometrically on each of these trees. We decompose the set $\{1,\ldots , k\}$ according to the properties of the action $\Gamma_{0} \curvearrowright T_{i}$. Write
$$\{1,\ldots , k\}:= I_{1} \cup I_{2} \cup I_{3}$$
where $I_{1}$ is the set of indices $i$ such that $\Gamma_{0}$ fixes a point on $T_{i}$, $I_{2}$ is the set of indices $i$ such that $\Gamma_{0}$ preserves a finite set in the boundary $\partial T_{i}$ of $T_{i}$ but no point in $T_{i}$ itself, finally $I_{3}$ is the set of remaining indices.

\begin{lemma} For each $i\in I_{2}$, there exists a finite index subgroup $\Gamma_{i}\subset \Gamma_{0}$ and a homomorphism $\phi_{i} : \Gamma_{i}\to \R$ with the following property: each element in the kernel $H_{i}$ of $\phi_{i}$ fixes a point in $T_{i}$. 
\end{lemma}
{\it Proof of the lemma.} Let $F\subset \partial T_{i}$ be a finite $\Gamma_{0}$-invariant subset. A finite index subgroup $\Gamma_{i}$ of $\Gamma_{0}$ fixes $F$ pointwise. Let $b_{\xi} : \Gamma_{i}\to \R$ be the Busemann character associated to any point $\xi \in F$. Its kernel is made up of elements acting as elliptic isometries on $T_{i}$.\hfill $\Box$

We define:
$$\Gamma_{1}=\bigcap_{i\in I_{2}}\Gamma_{i}.$$
\noindent This group has finite index in $\Gamma$. We now deal with the actions on the trees $T_{i}$ for $i\in I_{3}$. In the following, $X_{1}$ is the finite cover of the K\"ahler manifold $X$ with fundamental group $\Gamma_{1}$.

\begin{prop}\label{cne} For each $i\in I_{3}$, there exists a fibration $X_{1}\to \Sigma_{i}$ such that the kernel $H_{i}$ of the induced map $(p_{i})_{\ast}: \Gamma_{1}\to \pi_{1}^{orb}(\Sigma_{i})$ fixes a point in $T_{i}$. 
\end{prop}
{\it Proof of the proposition.} This result is due to Gromov and Schoen~\cite{gs}. Here we only sketch the ideas of the proof, see~\cite[\S 6.6]{abckt} and~\cite{simpson} for details. Since the action of $\Gamma_{1}$ on $T_{i}$ is non-elementary (i.e. does not preserve any finite set in $T_{i}\cup \partial T_{i}$), there exists an equivariant pluriharmonic map $f : \widetilde{X}_{1}\to T_{i}$, where $\widetilde{X}_{1}$ is the universal cover of $X_{1}$ (see~\cite{gs}). This map gives rise to a (singular) holomorphic foliation of codimension $1$ on $X_{1}$ and one proves that this foliation is induced by a holomorphic fibration $p_{i}$ onto some hyperbolic $2$-dimensional orbifold $\Sigma_{i}$. The harmonic map $f$ is constant on the fibers of $p_{i}$, hence the kernel $H_{i}$ of the map

$$\Gamma_{1}\to \pi_{1}^{orb}(\Sigma_{i})$$

\noindent fixes pointwise the image of $f$ in $T_{i}$. This proves the proposition.\hfill $\Box$

\begin{rem} The trees constructed from the Davis complex of $W$ need not be locally finite. But this does not affect the proof of the previous proposition. As suggested to us by Marc Burger, one can also recover the result of Gromov and Schoen describing non-elementary actions of K\"ahler groups on trees as follows. One combines the fact that an action on a tree gives rise to an action on the infinite dimensional real hyperbolic space $\mathbb{H}^{\infty}_{\mathbb{R}}$~\cite{bim} with the description of actions of K\" ahler groups on the space $\mathbb{H}^{\infty}_{\mathbb{R}}$ obtained in~\cite{delzantpy}.  
\end{rem}

We now define:

$$H=\bigcap_{i\in I_{2}\cup I_{3}} H_{i}\subset \Gamma_{1}.$$
Each element of the group $H$ fixes a point on each of the trees $(T_{i})_{1\le i \le k}$: for $i\in I_{1}$ this is because $\Gamma_{1}$ itself fixes a point on $T_{i}$, for $i\in I_{2}\cup I_{3}$, this follows from the definition of the groups $H_{i}$. Since the action of $W_{0}$ on $T_{1}\times \cdots \times T_{k}$ is free, the group $\varphi (H)\subset W_{0}$ must be trivial. In other words, the restriction of $\varphi$ to $\Gamma_{1}$ factors through the homomorphism
$$\Gamma_{1} \to (\Gamma_{1})_{{\rm ab}}\times \prod_{i\in I_{3}} \pi_{1}^{orb}(\Sigma_{i}).$$ This concludes the proof of Theorem~\ref{kc}.\hfill $\Box$

Before proving Corollaries~\ref{coxcom} and~\ref{corraag}, let us recall the main result from~\cite{bhms2002}:

\begin{center}
{\it Let $G$ be a subgroup of a direct product $H_{1}\times \cdots \times H_{n}$ where each $H_{i}$ is a surface group or a free group. If $G$ is of type ${\rm FP}_{n}$, it has a subgroup of finite index isomorphic to a direct product of the form $A_{1}\times \cdots \times A_{r}$ where $r\le n$ and each $A_{i}$ is a finitely generated subgroup of one of the $H_{j}$'s.}
\end{center} 

Recall that a group $G$ is of type ${\rm FP}_{n}$ if there is an exact sequence 
$$P_{n}\to P_{n-1}\to \cdots \to P_{0}\to \Z \to 0$$
\noindent of $\Z G$-modules, where the $P_{i}$ are finitely generated and projective and where $\Z$ is considered as a trivial $\Z G$-module. See~\cite[\S VIII.5]{browncoho} for more details on this notion. We will apply this result to torsion-free finite index subgroups of Coxeter groups (since they act cocompactly and freely on the Davis complex, they admit a classifying space which is a finite complex, hence are of type ${\rm FP}_{\infty}$). For examples of K\"ahler groups which are not of type ${\rm FP}_{\infty}$, see~\cite{dps}.

Note that the result of~\cite{bhms2002} applies in particular to subgroups of direct products of the form
$$\mathbb{Z}^{l}\times H_{1} \times \cdots \times H_{m}$$
where the $H_{i}$'s are surface groups.

{\it Proof of Corollary~\ref{coxcom}.} Let $W$ be a Coxeter group admitting a finite index subgroup $H$ isomorphic to a K\"ahler group. According to Theorem~\ref{kc}, a finite index subgroup $H_{1}$ of $H$ admits a faithful homomorphism
$$\phi : H_{1}\to \Z^{l}\times \pi_{1}(\Sigma_{g_{1}})\times \cdots \times \pi_{1}(\Sigma_{g_{m}}),$$

\noindent where the $\Sigma_{g_{j}}$ are closed orientable surfaces of genus greater than $1$. Let $W_{i}$ be an infinite irreducible factor of $W$ which is not Euclidean. We will show that $W_{i}$ has a finite index subgroup isomorphic to a surface group. We start with the following lemma.

\begin{lemma}\label{free} The group $W_{i}$ is not virtually free. 
\end{lemma}
{\it Proof of the lemma.} Assume by contradiction that a finite index subgroup of $W_{i}$ is free of rank $\ge 2$ (note that, being non-Euclidean, $W_{i}$ is not virtually Abelian; this follows for instance from~\cite{bedeha}). We know that the group $W$ has a finite index subgroup $H$ which is a K\"ahler group. There is a finite index subgroup $H_{2}$ of $H$ which is a direct product of finite index subgroups of each irreducible factor of $W$. Hence, under our hypothesis, we can take $H_{2}$ of the form $F \times A$ where $F$ is free non-Abelian. But there is no K\"ahler group of the form $F\times A$ according to~\cite[Theorem 3]{jr1987}.\hfill $\Box$

Let $G:=H_{1}\cap W_{i}$.  The restriction of $\phi$ to $G$ gives a faithful homomorphism 
$$G\to \Z^{l}\times \pi_{1}(\Sigma_{g_{1}})\times \cdots \times \pi_{1}(\Sigma_{g_{m}}).$$

\noindent According to the result from~\cite{bhms2002} stated above, we obtain that a finite index subgroup $G_{1}$ of $G$ is isomorphic to a product 
$$A_{0}\times A_{1} \times \cdots \times A_{r}$$
where $r\le m$, $A_{0}$ is free Abelian and each $A_{i}$ ($1\le i \le r$) is a finitely generated subgroup of one of the $\pi_{1}(\Sigma_{g_{j}})$. By Proposition~\ref{nn}, there is only one nontrivial factor in this decomposition. Since $W_{i}$ is not virtually Abelian, this implies that $G_{1}$ is isomorphic to a subgroup of one of the $\pi_{1}(\Sigma_{g_{j}})$. Since $G_{1}$ cannot be free, according to Lemma~\ref{free}, it has to be of finite index in $\pi_{1}(\Sigma_{g_{j}})$. This proves the corollary.\hfill $\Box$

In the proof of the next corollary, we will use several times the following fact: if $A(\mathscr{G})$ is a RAAG and if $\mathscr{G}_{1}\subset \mathscr{G}$ is the subgraph with vertex set $V_{1}$, the subgroup of $A(\mathscr{G})$ generated by the $g_{v}$'s for $v\in V_{1}$ is isomorphic to the RAAG $A(\mathscr{G}_{1})$ (see~\cite[\S 3.2]{charney}). In particular, a pair of generators generates either a free group or a free Abelian group. 

{\it Proof of Corollary~\ref{corraag}.} Any RAAG $A(\mathscr{G})$ has a natural quotient which is a right-angled Coxeter group $W(\mathscr{G})$: one simply adds the relations $g_{v}^{2}=1$ to the presentation of the group. We will say that a RAAG is irreducible if $W(\mathscr{G})$ is irreducible. Any RAAG $A(\mathscr{G})$ can be written as a direct product of irreducible RAAGs
$$A(\mathscr{G})\simeq A(\mathscr{G}_{1})\times \cdots \times A(\mathscr{G}_{r}),$$

\noindent see~\cite{charneydavis}, Lemma~2.2.6. According to~\cite{dj-jpaa}, each irreducible factor $A(\mathscr{G}_{j})$ embeds as a subgroup of finite index in a Coxeter group $W_{j}$. One sees from the proof in~\cite{dj-jpaa} that the irreducibility of $A(\mathscr{G}_{j})$ implies that $W_{j}$ is also irreducible. 

If $A(\mathscr{G})$ is commensurable with a K\"ahler group, the Coxeter group $W_{1}\times \cdots \times W_{r}$ is also commensurable with a K\" ahler group. According to Corollary~\ref{coxcom}, each group $A(\mathscr{G}_{j})$ must be either virtually Abelian or virtually a surface group. The surface group case does not occur (if a RAAG does not contain $\Z^{2}$, it is free). Each factor $A(\mathscr{G}_{j})$ is thus virtually Abelian hence Abelian. This proves that $A(\mathscr{G})$ is Abelian, its rank being necessarily even if it is commensurable with a K\"ahler group.\hfill $\Box$

\begin{rem}\label{rzd} To prove Corollaries~\ref{coxcom} and~\ref{corraag}, one can also use the following argument in replacement of the result from~\cite{bhms2002}: if a group $G$ admits a Zariski dense embedding into a simple Lie group with trivial center, any two nontrivial normal subgroups of $G$ have nontrivial intersection. As a consequence, if $G$ embeds into a direct product $A\times B$, $G$ embeds in either $A$ or $B$. This applies to irreducible, infinite, non-Euclidean Coxeter groups (as follows from the results in~\cite{bedeha} and~\cite[\S 6.1]{krammer}). This alternative proof was pointed out to us by Yves de Cornulier.  
\end{rem}


\bigskip
\bigskip
\begin{small}
\begin{tabular}{l}
IRMA, Universit\'e de Strasbourg \& CNRS\\
67084 Strasbourg, France\\
ppy@math.unistra.fr\\    
\end{tabular}
\end{small}

 \end{document}